\documentclass[reqno, a4paper, 10pt,oneside]{amsart}

\makeatletter
\def\specialsection{\@startsection{section}{1}%
	\z@{\linespacing\@plus\linespacing}{.5\linespacing}%
	{\normalfont}}
\def\section{\@startsection{section}{1}%
	\z@{.7\linespacing\@plus\linespacing}{.5\linespacing}%
	{\normalfont\scshape\bfseries}}

\makeatother

\usepackage[top=1in, bottom=1.4in, inner=1in, outer=1in, includehead, showframe=false]{geometry}
\usepackage[algo2e, linesnumbered, noline, noend, ruled]{algorithm2e}
\usepackage[foot]{amsaddr}
\usepackage{amsfonts}
\usepackage{amsmath}
\usepackage{amssymb}
\usepackage{amsthm}
\usepackage[english]{babel}
\usepackage{booktabs}
\usepackage{caption}
\usepackage{braket}
\usepackage{float}
\usepackage[T1]{fontenc}
\usepackage{graphicx}
\usepackage[utf8]{inputenc}
\usepackage{lmodern}
\usepackage{lipsum}
\usepackage{mathtools}
\usepackage{nicefrac}
\usepackage[defaultlines=2,all]{nowidow}
\usepackage[per-mode=reciprocal-positive-first,group-digits=integer,range-units=single]{siunitx}
\usepackage[table,dvipsnames]{xcolor}
\usepackage{enumitem}
\usepackage{subcaption}
\usepackage{stmaryrd}
\usepackage[version=4]{mhchem}
\usepackage{textcomp}

\usepackage{cite}
\usepackage{tikz}
\usepackage[allcolors=cyan,colorlinks]{hyperref}
\usepackage[english,capitalize,noabbrev]{cleveref}

\DeclareSIUnit\bar{bar}
\DeclareSIUnit{\molepercent}{mol\%}

\definecolor{myblue}{RGB}{0, 107, 164}
\definecolor{myorange}{RGB}{255, 128, 14}
\definecolor{grayone}{RGB}{171, 171, 171}

	\theoremstyle{plain}
	\newtheorem{Theorem}{Theorem}[section]

	\theoremstyle{definition}
	
	\newtheorem{Remark}[Theorem]{Remark}
	
	\newtheorem*{Example*}{Example}
	\newtheorem*{Remark*}{Remark}
	
	\setenumerate{label=(\arabic*), ref=(\arabic*)}
	
	\makeatletter
	\DeclareOldFontCommand{\rm}{\normalfont\rmfamily}{\mathrm}
	\DeclareOldFontCommand{\sf}{\normalfont\sffamily}{\mathsf}
	\DeclareOldFontCommand{\tt}{\normalfont\ttfamily}{\mathtt}
	\DeclareOldFontCommand{\bf}{\normalfont\bfseries}{\mathbf}
	\DeclareOldFontCommand{\it}{\normalfont\itshape}{\mathit}
	\DeclareOldFontCommand{\sl}{\normalfont\slshape}{\@nomath\sl}
	\DeclareOldFontCommand{\sc}{\normalfont\scshape}{\@nomath\sc}
	\makeatother

	\numberwithin{equation}{section}
	
	\newcommand*\closure[1]{%
		\hbox{%
			\vbox{%
				\hrule height 0.5pt 
				\kern0.5ex
				\hbox{%
					\ensuremath{#1}%
				}%
			}%
		}%
	} 
	
\SetCommentSty{mycommfont}
\SetKwComment{Comment}{\# }{}

\makeatletter
\ifx\insert\@title\relax
  \ClassWarning{\@classname}{You need to use \string\title\space to set the title.}
\fi
\renewenvironment{abstract}{%
  \ifx\maketitle\relax
    \ClassWarning{\@classname}{Abstract should precede \maketitle.}%
  \fi
  \global\setbox\abstractbox=\vbox\bgroup
  \vspace{1em}
  \normalfont\Small
  \list{}{%
    \labelwidth\z@ \leftmargin3pc \rightmargin\leftmargin
    \listparindent\normalparindent \itemindent\z@
    \parsep\z@ \@plus\p@
    
  }%
  \item[]%
}{%
  \endlist\egroup
  \ifx\@setabstract\relax \@setabstracta \fi
}
\makeatother

\begin{document}

\title[CFD-Based Shape Optimization of Structured Packings for Distillation]{CFD-based Shape Optimization of Structured Packings for Enhancing Separation Efficiency in Distillation}
\author{Sebastian Blauth$^{*,1}$}
\author{Dennis Stucke$^{2}$}
\author{Mohamed Adel Ashour$^{2}$}
\author{Johannes Schnebele$^{1}$}
\author{Thomas Gr\"utzner$^{2}$}
\author{Christian Leith\"auser$^1$}

\address{$^*$ Corresponding Author}
\address{$^1$ Fraunhofer ITWM, Department Transport Processes, Kaiserslautern, Germany}
\address{$^2$ Ulm University, Institute of Chemical Engineering, Laboratory of Thermal Process Engineering, Ulm, Germany}
\email{\href{mailto:sebastian.blauth@itwm.fraunhofer.de}{sebastian.blauth@itwm.fraunhofer.de}}

\begin{abstract}
    Free-form shape optimization techniques are investigated to improve
    the separation efficiency of structured packings in laboratory-scale distillation columns. A simplified simulation model based on computational fluid dynamics (CFD) for the mass transfer in the distillation column is used and a corresponding shape optimization problem is formulated. The goal of the optimization is to increase the mass transfer in the column by changing the packing's shape, which has been previously used as criterion for increasing the separation efficiency of the column. The computational shape optimization yields promising results, with an increased mass transfer of nearly \qty{20}{\percent}. For validation, the resulting optimized shape is additively manufactured using 3D-printing and investigated experimentally. The experimental results are in good agreement with the performance improvement predicted by the computational model, yielding an increase in separation efficiency of around \qty{20}{\percent}. 
 
	\medskip
	\noindent \textsc{Keywords. } Shape Optimization, Additive Manufacturing, Structured Packing, Distillation, CFD
	
\end{abstract}


{\noindent\footnotesize This is a post-peer-review, pre-copyedit version of an article published in Chemical Engineering Science. The final version is available online at \url{https://doi.org/10.1016/j.ces.2024.120803}.
}

\maketitle

\vspace{-0.8cm}

\section{Introduction}
\label{sec:introduction}

In past years the research in the field of structured packing development for laboratory-scale separation processes has intensified. In many cases 3D-printing is used to manufacture these new structures.~\cite{neukaufer_development_2022,al-maqaleh_experimental_2022,dejean_design_2020,flagiello_performances_2023,sun_hydrodynamics_2021,olenberg_optimization_2018} The main objective is to miniaturize laboratory columns regarding the column diameter. This reduction in diameter has several advantages like the decreased size of the column itself and the peripherals. Another important benefit is the reduced cost for chemicals due to the lower holdup in the column and thereby decreased operational costs. This also benefits the time-to-market as there is less time needed to synthesize the required chemicals in laboratory experiments. Additionally, the reduced amount of chemicals lowers the safety requirement for laboratory experiments. \cite{neukaufer_flexible_2021,norris_precise_1945} On the other hand this reduction in diameter causes additional problems related to the increased surface-to-volume ratio of the cylindrical column. One negative effect is the stronger impact of heat losses on the internal flow rate in the column, which has been investigated by Ashour et al.~\cite{ashour_revisiting_2023}. The other complication is the liquid maldistribution which leads to increased liquid flow rates close to the wall and reduced liquid flow rates in the middle of the packing. To overcome this issue, the packing structure can be designed in a way to counteract this effect \cite{neukaufer_development_2022}.

In the literature researches present different approaches to design structured packings. Neuk\"{a}ufer et al.~\cite{neukaufer_development_2022} designed a repetitive unit cell by using conventional computer-aided-design (CAD) programs. This unit cell is then repeated in all directions and, finally, cut into cylindrical shape. Miramontes et al.~\cite{miramontes_additively_2020} used a similar approach, but instead of starting with small unit cells they designed structured sheets which are then stacked and cut to obtain the cylindrical shape. Zimmer et al.~\cite{zimmer_effect_2021} used triple periodic minimal surface to create packing structures. G\l{}adyszewski et al.~\cite{gladyszewski_additive_2018} used micro computer tomography to scan an existing metallic foam and recreate it by additive manufacturing. All these methods have in common that the development of new structures and the improvement of existing structures is based on educated guesses by the engineers.

In this article a novel approach to enhance the separation efficiency is presented. For this, a CAD-model of an existing structure is required as starting point. The selected structure is the 3D-printable version of the Rombopak 9M (RP9M-3D) which is presented in~\cite{neukaufer_development_2022}. The shape of this packing is then optimized using free-form shape optimization based on a single-phase CFD surrogate model for the distillation column.

The shape optimization approach considered here is a CAD-free or free-form shape optimization. This means that the shape is not parametrized in some CAD model and then the parameters of this model are optimized, but instead the nodes of the computational mesh are moved to alter the shape. In contrast to the parametrized approach, CAD-free shape optimization does not restrict the reachable shapes and is, thus, more flexible. For further information, the reader is referred, e.g., to \cite{Delfour2011Shapes, Sokolowski1992Introduction}.

Shape optimization has been used to improve designs for many applications. These include, e.g., the optimization of aircraft \cite{Schmidt2013Three}, electrochemical cells \cite{Blauth2024Multi} or microchannel systems \cite{Blauth2021Model}. Recently, shape optimization has received a lot of attention, e.g., in \cite{Mueller2021novel}, where special mesh deformations for shape optimization are investigated, in \cite{Schulz2016Efficient,Blauth2021Nonlinear}, where efficient shape optimization algorithms are proposed, or in \cite{Blauth2023Space}, where space mapping methods for shape optimization are presented. In this paper, it is investigated how these techniques can be applied to optimize the shape of packings for distillation columns. To the best of our knowledge, using techniques from shape optimization to improve the design of structured packings for distillation has not yet been considered in the literature.

\section{Materials and Methods}
\label{sec:methods}

\subsection{Packing structure}
\label{ssec:packing}

This article focuses on presenting a shape optimization approach for packing structures and validating it with experimental methods. A Sulzer Rombopak 9M, which was first modified to be manufacturable by additive manufacturing \cite{neukaufer_development_2022}, is used as base structure. The modification consists of increasing the thickness of all parts of the structure to \qty{0,8}{\milli\meter} to meet the 3D-printer's minimum requirements. This is possible due to the large void fraction of the chosen structure. The resulting structure is called RP9M-3D and in \cref{fig:rp9m_cad} a CAD image of this structure is presented.

There are multiple reasons for choosing the RP9M-3D as base structure. Most packings that are used in laboratory-scale experiments for scale-up purposes are sheet metal packings like the Rombopak structure. These packings generally feature a rather constant separation efficiency over the operation range which is a crucial property for packings used to generate scale-up data. The second objective for packings for scale-up purposes is the increase of separation efficiency, but without obstructing the main objective mentioned before. While there are other types of packings, such as wire mesh ones, that show a higher separation efficiency than the considered sheet metal ones, these packings can usually not be used for scale up purposes due to their non-constant behavior over their operation range.

\subsection{Experimental Methods}
\label{ssec:experimental_methods}

\begin{figure}[t]
    \centering
    \begin{subfigure}[b]{0.4\textwidth}
        \centering
        \includegraphics[angle=-90, origin=c, height=0.25\textheight]{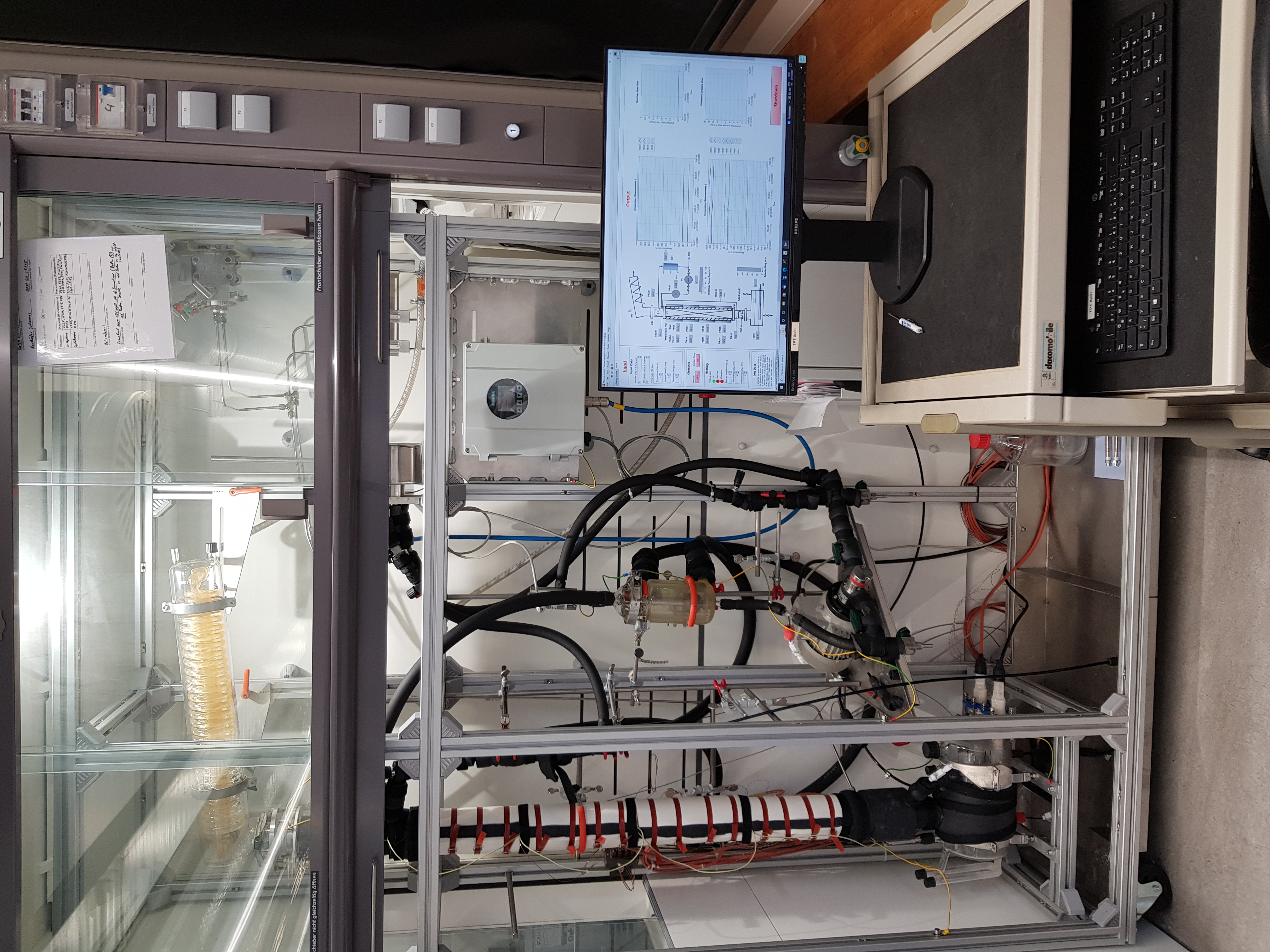}
        \caption{Laboratory-scale batch distillation column used for the measurement of the separation efficiency of structured packings.}
        \label{fig:test-rig}
    \end{subfigure}%
    \hfill
    \begin{subfigure}[b]{0.3\textwidth}
        \centering
        \includegraphics[height=0.3\textheight]{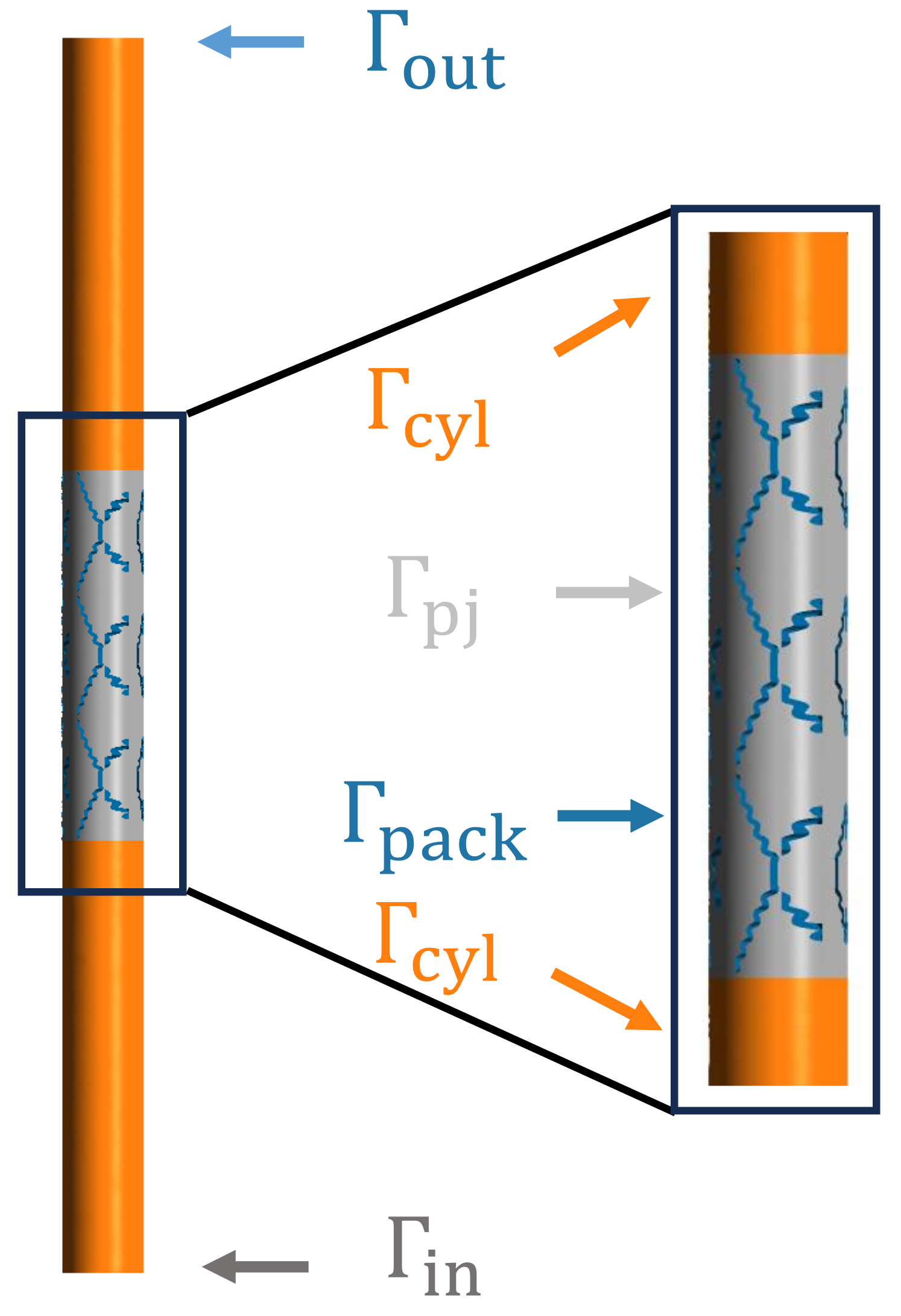}
        \caption{Geometrical setup of the model and boundary conditions.}
        \label{fig:geometry}
    \end{subfigure}%
    \hfill
    \begin{subfigure}[b]{0.25\textwidth}
        \centering
        \includegraphics[height=0.3\textheight]{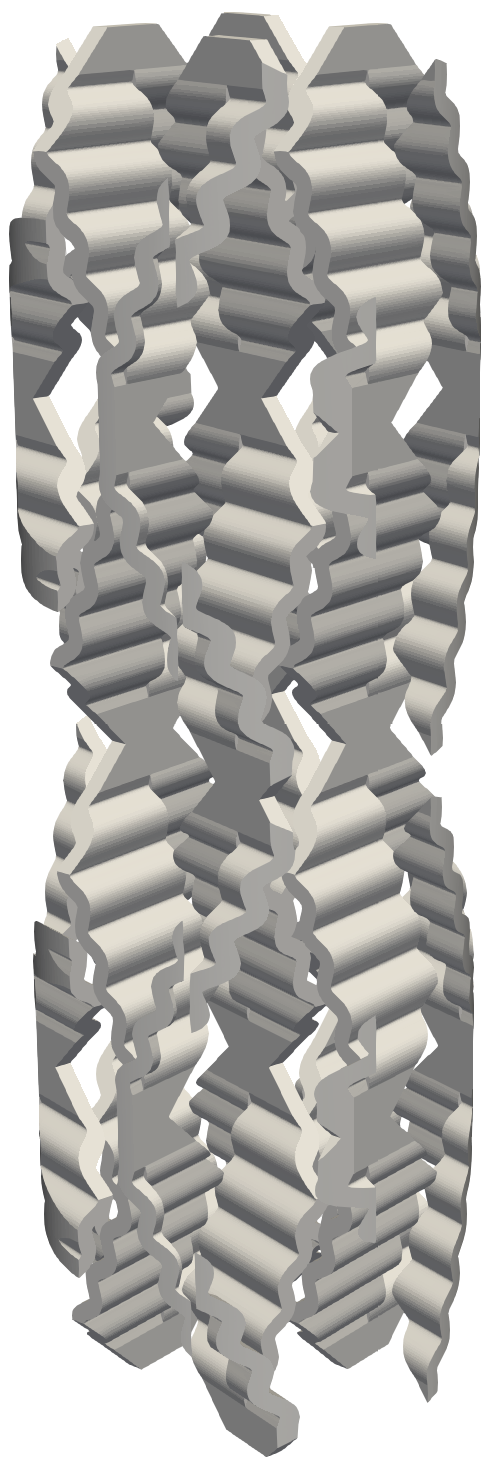}
        \caption{CAD image showing the RP9M-3D packing structure.}
        \label{fig:rp9m_cad}
    \end{subfigure}
    \caption{Experimental and simulation setup.}
    \label{fig:setup}
\end{figure}

The experimental investigation of the separation efficiency is done in a laboratory-scale distillation column which is shown in \cref{fig:test-rig}. As this test-rig is explained in detail in previous publications \cite{neukaufer_flexible_2021,ashour_revisiting_2023}, only a brief description is given here. The device is designed to precisely measure the separation efficiency and heat losses of different structured packings at various operating conditions. Therefore, it is operated at batch mode under total reflux condition. As test system the standard system of cyclohexane/\textit{n}-heptane\cite{onken_recommended_1989} is used with a concentration of \qty{25}{\molepercent} cyclohexane. The column consists of an additively manufactured column section combined with custom designed multifunctional trays (MFT). The MFT's serve as phase separators for gas and liquid phase which allows for sampling and the measurement of the flow rate. The MFT's and the column sections are made out of polyamid 12 (PA12) and are manufactured with the selective laser sintering (SLS) three-dimensional (3D) printer Formlabs Fuse1. The auxiliary parts are standard laboratory equipment like an electrically heated glass reboiler with \qty{2}{\kilo\watt} heating duty, a three-coil glass total condenser and high precision sensors to measure the temperature, pressure drop, and flow rate. Additionally the test rig features an active insulation which was implemented to reduce the heat losses through the wall. The process control system is implemented in LabVIEW\texttrademark. The composition of the samples is measured using the gas chromatograph Shimadzu GC-2030 ATF with a semi polar column and a flame ionization detector. To keep the influence of the measuring error small, it is ensured that the measured concentrations are not in the pinch zones of the light boiler composition. A detailed explanation of the measuring procedure, an investigation of its repeatability, and the equations used to calculate the separation efficiency can be found in \cite{neukaufer_development_2022,neukaufer_flexible_2021,ashour_revisiting_2023}.

\subsection{Simulation Methods}
\label{ssec:simulation_methods}

\subsubsection{Simulation Model}
\label{ssec:simulation_model}

Although sophisticated two-phase computational fluid dynamics (CFD) simulations of distillation columns are available \cite{neukaufer_development_2022}, the computational complexity of such models is prohibitively high for employing them for optimizing the packing's shape. For these reasons, a simplified single-phase CFD model from \cite{neukaufer_development_2022} is employed as surrogate model for the optimization, which is briefly recalled in the following. The simplified model considers the steady-state incompressible Navier-Stokes equations for the gas-phase flow and uses a convection-diffusion equation for modeling the mass transfer over the packing. The former is given by
\begin{equation}
	\label{eq:navier_stokes}
	\begin{aligned}
		-\mu \Delta u + \rho (u\cdot \nabla) u + \nabla p &= 0 \quad &&\text{ in } \Omega, \\
		\nabla \cdot u &= 0 \quad &&\text{ in } \Omega, \\
		%
		%
		%
	\end{aligned}
\end{equation}
with the boundary conditions
\begin{equation*}
    u = u_\mathrm{in} \quad \text{ on } \Gamma_\mathrm{in}, \qquad\quad u = 0 \quad \text{ on } \Gamma_\mathrm{wall}, \qquad\quad \mu \partial_n u - p n = 0 \quad \text{ on } \Gamma_\mathrm{out},
\end{equation*}
where $u$ and $p$ denote the fluid velocity and pressure, $\mu$ is the dynamic viscosity, and $\rho$ the density. The fluid domain is denoted by $\Omega \subset \mathbb{R}^3$ and its boundary $\Gamma$ is divided into the inlet $\Gamma_\mathrm{in}$, where the velocity is specified to be $u_\mathrm{in}$, the wall boundary $\Gamma_\mathrm{wall}$, where a no-slip condition holds, and the outlet $\Gamma_\mathrm{out}$, where the flow satisfies a do-nothing condition \cite{Elman2014Finite}.  The wall boundary is divided into the cylinder jacket $\Gamma_\mathrm{cyl}$, the packing jacket $\Gamma_\mathrm{pj}$, and the packing $\Gamma_\mathrm{pack}$ itself. Additionally, $n$ is the outer unit normal vector on $\Gamma$ and $\partial_n u$ is the normal derivative of $u$. The geometrical setup is shown in Figure~\ref{fig:geometry}. 

To model the mass transfer, the following convection-diffusion equation is used
\begin{equation}
	\label{eq:concentration}
	\begin{aligned}
		-\nabla\cdot (D \nabla c) + u \cdot \nabla c &= 0 \quad &&\text{ in } \Omega, \\
		%
		%
		%
	\end{aligned}
\end{equation}
with the boundary conditions
\begin{equation*}
    c = c_\mathrm{in} \quad \text{ on } \Gamma_\mathrm{in}, \qquad\quad c = c_\mathrm{pack} \quad \text{ on } \Gamma_\mathrm{pj} \cup \Gamma_\mathrm{pack}, \qquad\quad D \partial_n c = 0 \quad \text{ on } \Gamma_\mathrm{cyl} \cup \Gamma_\mathrm{out},
\end{equation*}
where $c$ is a fictitious concentration and $D$ the diffusion coefficient. Note that the concentration is fictitious in the sense that it does not resemble a physical quantity but is solely used for characterizing the mass transfer in the surrogate CFD model. The (arbitrary) unit of \unit{\mol\per\cubic\meter} is chosen for the concentration. For a physical interpretation of the simplified CFD model the reader is referred to Remark~\ref{rem:simplified_model} below.
The material properties are based on those of nitrogen at \qty{25}{\degreeCelsius} and \qty{1}{\bar}, i.e., $\rho~=~\qty{1.138}{\kilogram\per\cubic\meter}$ and $\mu~=~\qty{1.728e-5}{\pascal\second}$ \cite{stephan_vdi-warmeatlas_2019}. For the diffusion coefficient, the gas-side diffusion coefficient of the equimolar cyclohexane/\textit{n}-heptane system at \qty{1}{\bar} and boiling temperature of \qty{90}{\degreeCelsius}, i.e., $D = \qty{3.72e-6}{\square\meter\per\second}$, is used \cite{neukaufer_development_2022}. For the boundary conditions, values of $u_\mathrm{in} = \qty{0.933}{\meter\per\second}$, $c_\mathrm{in} = \qty{100}{\mol\per\cubic\meter}$, and $c_\mathrm{pack}~=~\qty{1}{\mol\per\cubic\meter}$ are used. 

To investigate the mass transfer, the logarithmic mass transfer coefficient $\beta$ \cite{neukaufer_development_2022} is used:
\begin{equation}
    \label{eq:mtc}
	\beta = \frac{\dot{V}}{A_\mathrm{geo}} \log\left( \frac{c_\mathrm{pack} - c_\mathrm{in}}{c_\mathrm{pack} - c_\mathrm{out}} \right),
\end{equation}
where $c_\mathrm{out}$ is the flow-averaged outlet concentration, $A_\mathrm{geo}$ denotes the surface area of the packing and packing wall, and $\dot{V}$ is the volume flow rate of the fluid.
In our previous work \cite{neukaufer_development_2022}, the logarithmic mass transfer coefficient has been used as optimization criterion for improving the packing design and it has been found that the improvement predicted with this measure agrees well with experimental measurements of the column's separation efficiency. 

To establish a CFD model with which the results of the shape optimization can be validated numerically, Ansys\textsuperscript{\textregistered} Fluent \cite{Ansys2023Ansys} is employed for simulating \eqref{eq:navier_stokes} and \eqref{eq:concentration} with an SST $k-\omega$ turbulence model. Note that this slightly increases the diffusion coefficient $D$ due to the turbulent diffusivity \cite{neukaufer_development_2022}. Fluent is based on a finite volume method and used to recreate the simulation results from \cite{neukaufer_development_2022}, which were originally obtained with OpenFOAM\textsuperscript{\textregistered}. For a detailed description of the computational mesh and the solver settings for the simulation with Ansys\textsuperscript{\textregistered} Fluent the reader is referred to Appendix~\ref{app:fluent_settings}.

\begin{Remark}
    \label{rem:simplified_model}
    While the simplified CFD model presented here has been successfully applied for the optimization of structured packings in the previous literature \cite{neukaufer_development_2022}, it shall be emphasized that it is based on a series of simplifying assumptions. As described in \cite{neukaufer_development_2022}, the simplified model neglects the liquid-side mass transfer resistance, so that it can only give results that are qualitatively correct. Moreover, it is assumed that the entire packing surface as well as the packing jacket participate in the mass transfer and are, thus, wetted completely, which is not the case in reality. Additionally, the value of the fictitious concentration on the packing and packing jacket is constant over time. This resembles an infinite sink for the transferred species. Finally, the film thickness is assumed to be infinitesimally small so that the gas flow is not impeded. Again, this assumption is not valid in reality. However, note that under these assumptions, $\beta$ from \eqref{eq:mtc} can, in fact, be interpreted as the gas-phase (logarithmic) mass transfer coefficient in the simplified model of the distillation process.
    
    The simplified model cannot and is not intended to model the intricate interplay of liquid and gas phase present in distillation columns. It is merely used as computational tool which models some of the mass transfer aspects present in distillation processes and allows for a comparatively fast simulation. For this reason, it enables the shape optimization approach presented in the subsequent section. In particular, it shall be emphasized that it is currently not feasible to perform such a shape optimization using a high-fidelity two-phase simulation model of the distillation process, as the associated computational cost would be way too high. For these reasons, either using a high-fidelity two-phase CFD model or experiments to validate the shape optimization results is of utmost importance and is considered in Section~\ref{ssec:experimental_results}. Such a validation should be used to investigate the reliability of the optimization results in practice, particularly, to ensure that optimization approach does not introduce geometrical changes which are not considered by the simplified model (e.g.\ the liquid distribution) but could have a detrimental impact on the practical performance of the packing.
    
    It is subject of future research to refine the simplified computational model to consider more realistic effects of the two-phase flow and mass transfer. Additionally, it is planned to use two-phase CFD simulations to numerically evaluate designs which have been obtained by the presented shape optimization approach.
\end{Remark}

\subsubsection{Shape Optimization}
\label{ssec:shape_optimization}

To increase the separation efficiency of the packing, methods from shape optimization constrained by partial differential equations (PDEs) are employed. The corresponding shape optimization problem is given by
\begin{equation}
	\label{eq:shape_opt}
    \max_{\Omega, c} J(\Omega, c) = \beta(c) = \frac{\dot{V}}{A_\mathrm{geo}} \log\left( \frac{c_\mathrm{pack} - c_\mathrm{in}}{c_\mathrm{pack} - c_\mathrm{out}} \right) \quad \text{ subject to } \eqref{eq:navier_stokes} \text{ and } \eqref{eq:concentration},
	\begin{aligned}
	\end{aligned}
\end{equation}
where $\beta$ depends on the concentration $c$ via $c_\mathrm{out}$, cf.~\eqref{eq:mtc}. The aim of \eqref{eq:shape_opt} is to maximize $\beta$, which in turn should yield a better separation efficiency of the packing, as discussed previously. The following geometrical constraints for \eqref{eq:shape_opt} are used: First, the inlet $\Gamma_\mathrm{in}$, outlet $\Gamma_\mathrm{out}$, and the cylinder jacket $\Gamma_\mathrm{cyl}$ are fixed and cannot change their shape. Second, the packing jacket $\Gamma_\mathrm{pj}$ is fixed in the direction normal to the jacket, so that it can only change tangentially. These constraints ensure that only the packing surface $\Gamma_\mathrm{pack}$ is changed during the shape optimization. Note that a CAD-free shape optimization is considered, i.e., the geometry is not parametrized, but instead the location of the mesh nodes of the computational mesh is changed. For more details regarding shape optimization, the reader is referred, e.g., to \cite{Delfour2011Shapes}.

For the numerical solution of the shape optimization problem \eqref{eq:shape_opt} the software package cashocs \cite{Blauth2021cashocs,Blauth2023Version} is employed.
As cashocs is based on the finite element software FEniCS \cite{Alnes2015FEniCS}, equations \eqref{eq:navier_stokes} and \eqref{eq:concentration} also have to be simulated with FEniCS to facilitate the shape optimization with cashocs. To do so, piecewise linear Lagrange elements are used to discretize all variables. Moreover, PSPG, SUPG, and grad-div stabilization for \eqref{eq:navier_stokes} and SUPG and crosswind diffusion stabilization for \eqref{eq:concentration} are used (see, e.g., \mbox{\cite{Elman2014Finite,Roos1996Numerical}}). For more details regarding the discretization, stabilization, and numerical solution of equations \eqref{eq:navier_stokes} and \eqref{eq:concentration} the reader is referred to Appendix~\ref{app:fenics_settings}.

For the shape optimization, 90 steps of a gradient descent method (see, e.g., \cite{Blauth2021Nonlinear}) implemented in cashocs are considered together with mesh quality constraints \cite{Blauth2024Gradient} to ensure that the mesh quality does not deteriorate during the optimization. The initial shape is given by the RP9M-3D from \cite{neukaufer_development_2022}, see Figure~\ref{fig:setup}, and the computational mesh detailed in Appendix~\ref{app:fluent_settings}, which consists of about \num{13.7}~million tetrahedrons and \num{2.35}~million vertices, is used as initial mesh. The resulting shape optimization problem has about \num{7}~million optimization variables, namely the coordinates of the mesh's vertices, which makes it a very large-scale optimization problem that can only be treated efficiently with an adjoint approach, which is implemented in our software cashocs.
To validate the numerical results obtained with FEniCS for the shape optimization, the commercial solver Ansys\textsuperscript{\textregistered} Fluent is used as discussed in Section~\ref{ssec:simulation_model}. Hence, during the course of the shape optimization, the stabilized finite element formulation described above is used. After the optimization, both the initial and optimized designs are, again, simulated, this time using Ansys\textsuperscript{\textregistered} Fluent, to validate the numerical results and to ensure that the effects of the numerical stabilization in the finite element model is not too large.
As shape optimization problems have many local minimizers due to the non-convexity of the cost functional even for linear PDE constraints, the gradient descent method can only converge to a nearby local minimizer \cite{Allaire2021Shape}. For this reason, we expect that the optimized shape will not be too different from the initial one of the RP9M-3D.

\section{Results and Discussion}
\label{sec:results}

\subsection{Shape Optimization of the Packing}
\label{ssec:simulation_results}

\begin{figure}[!b]
	\centering
	\includegraphics[width=0.5\textwidth]{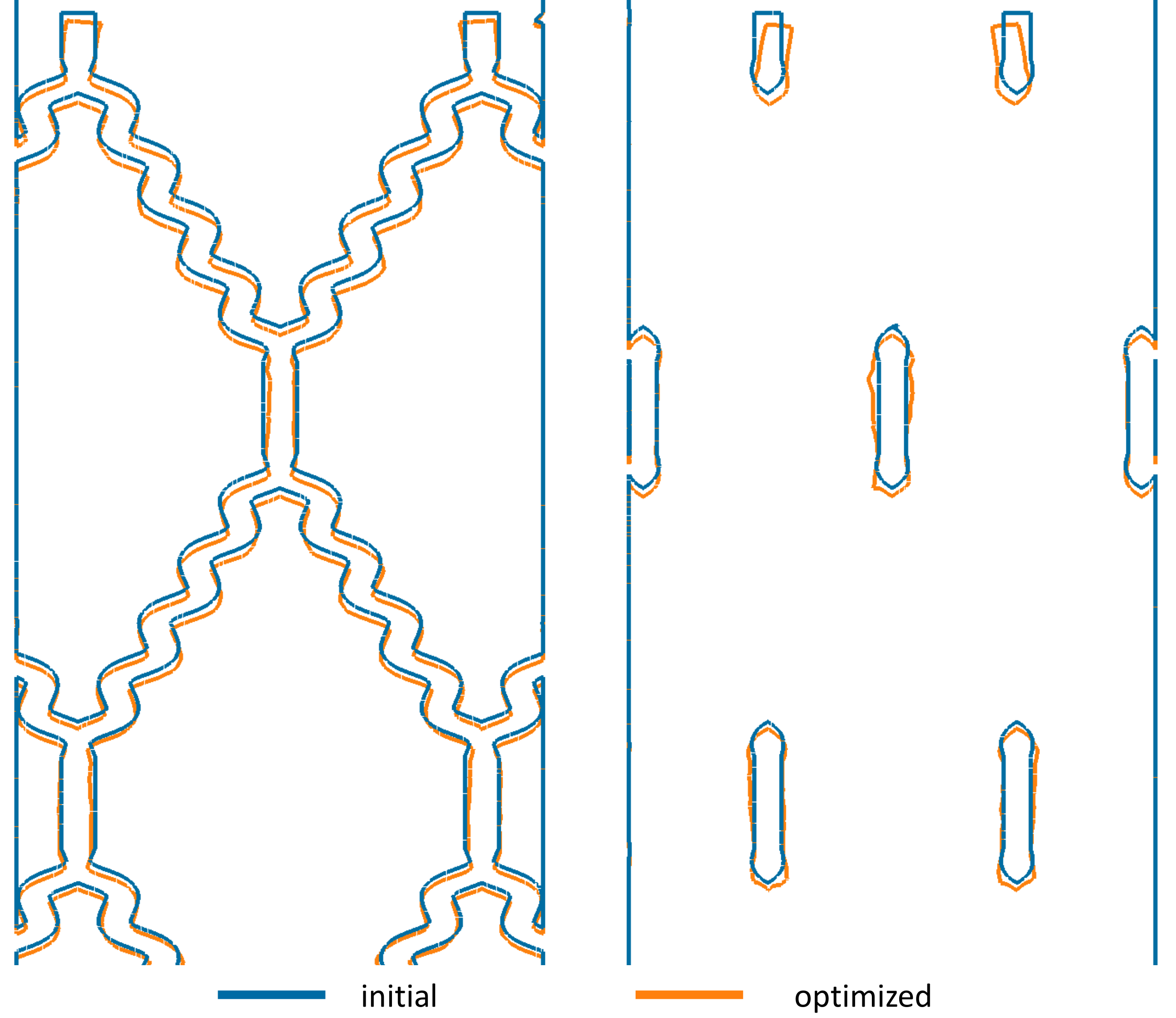}
	\caption{Comparison of the initial and optimized designs on the $xy$-plane (left) and $yz$-plane (right).}
	\label{fig:geometry_comparison}
\end{figure}

The shape optimization was run on 32 cores of an Intel\textsuperscript{\textregistered} Xeon\textsuperscript{\textregistered} Gold 6240R CPU and it took about 50~hours to complete the 90 iterations of the gradient descent method and find the optimized design. In Figure \ref{fig:geometry_comparison} the difference between initial and optimized designs is shown. The optimized shape does only change slightly in comparison to the initial one. E.g., it is observed that the stair-like structures in the $xy$-plane are slightly shifted downward and that the crosspieces in the $yz$-plane change their shape and thickness, particularly in the upwind direction. But, altogether, the changes between initial and optimized geometry are rather minute. This is attributed to the gradient descent method used for the shape optimization, which can, by design, only converge to a nearby local minimizer.

\begin{figure}[t]
    \centering
    \begin{subfigure}{0.475\textwidth}
        \centering
        \includegraphics[width=\textwidth]{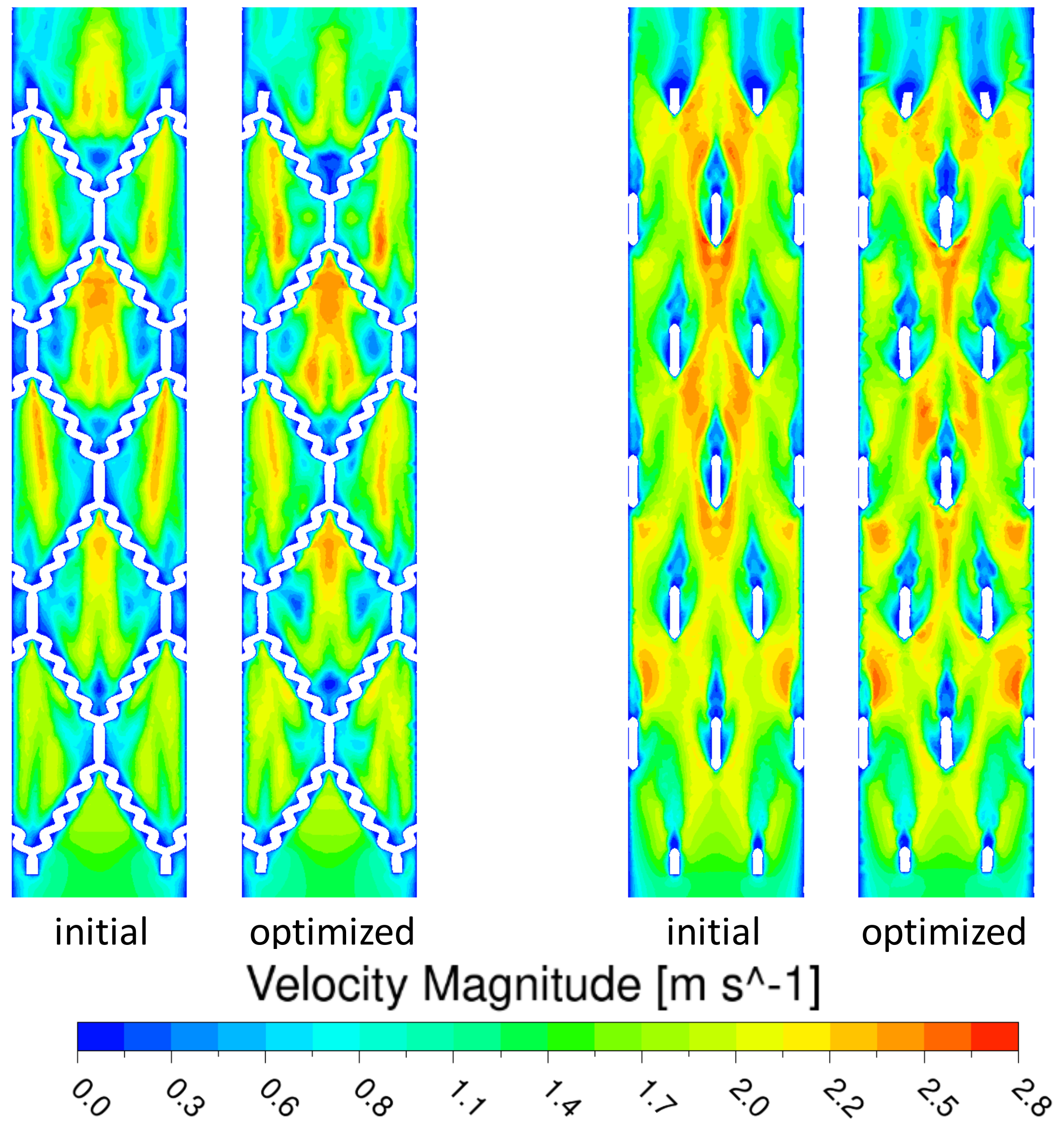}
        \caption{Velocity on initial and optimized designs on the $xy$- (left) and $yz$-plane (right).}
        \label{fig:velocity_comparison}
    \end{subfigure}
    \hfill
    \begin{subfigure}{0.475\textwidth}
        \centering
        \includegraphics[width=\textwidth]{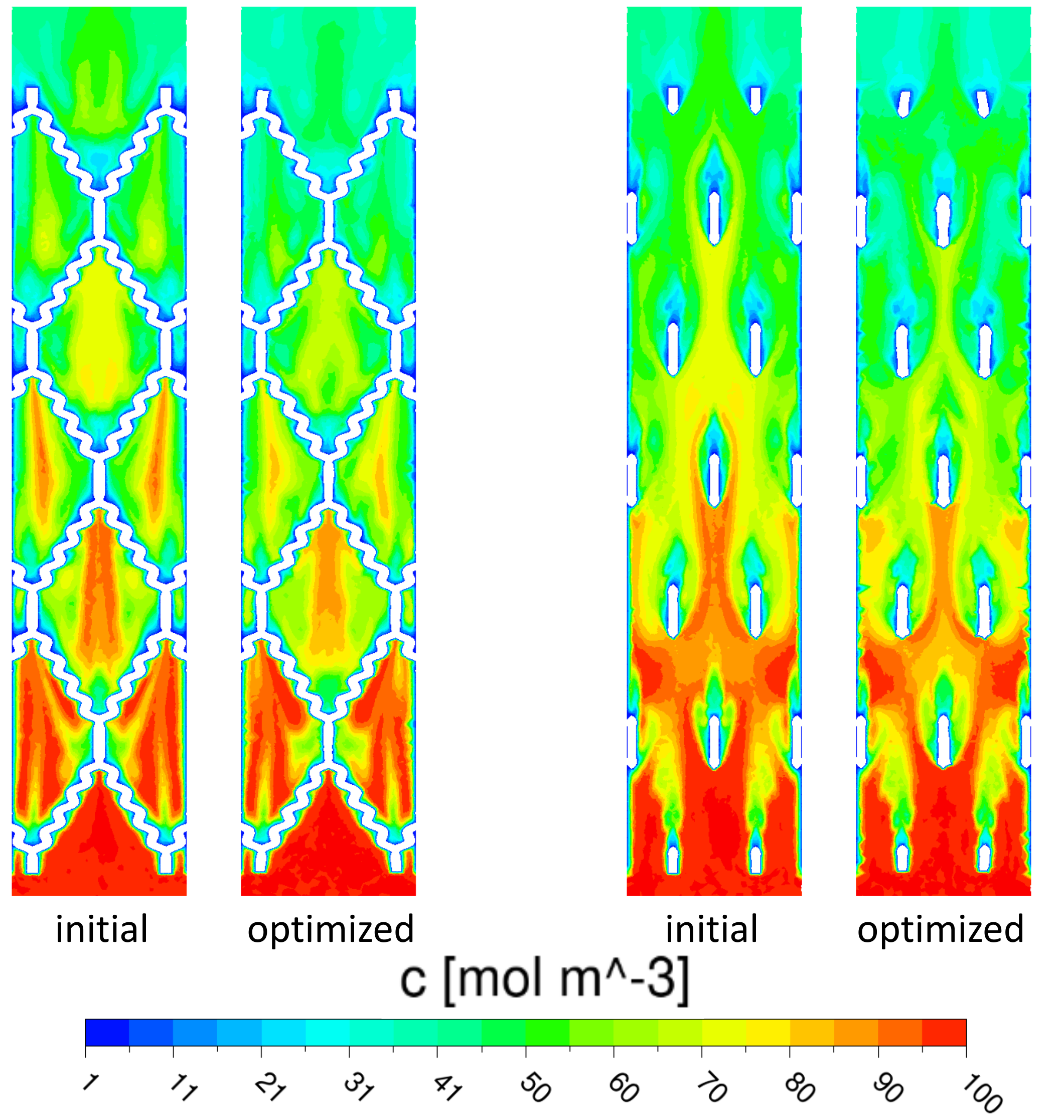}
        \caption{Concentration on initial and optimized designs on the $xy$- (left) and $yz$-plane (right).}
        \label{fig:concentration_comparison}
    \end{subfigure}
    \caption{Comparison of velocity (a) and pressure (b) on initial and optimized designs.}
    \label{fig:state_comparison}
\end{figure}

In Figure~\ref{fig:state_comparison}, the velocity and concentration are compared on the initial and optimized designs. It is observed that the small geometrical changes impact the velocity and concentration. 
In the $xy$-plane, the changes of the velocity magnitude are only minor. However, one can observe that the velocity magnitude is slightly increased locally, particularly near the stair-like structures of the packing wall, which results in an increase of the convective mass transfer near the wall. In the $yz$-plane, a different effect can be seen. There, it can be observed that the changes to the crosspieces lead to a deflection of the flow so that this is more spread out and results in the flow hitting more crosspieces downstream. This, too, leads to an increase of velocity near the packing surface and, consequently, to an increase of the convective mass transfer. The effects of these changes and the increased mass transfer can clearly be seen in the plots of the concentration, where we observe a significant decrease of the latter, indicating that the optimized packing has better mass transfer characteristics.
This is confirmed by an investigation of the logarithmic mass transfer coefficient: Compared to the initial value of $\beta=\qty{1.54e-2}{\meter\per\second}$, the optimized packing shows a substantial increase of \qty{19.7}{\percent} with $\beta=\qty{1.84e-2}{\meter\per\second}$. This change in $\beta$ indicates an enhanced mass transfer, particularly, as the geometric surface area of the packing only increased by \qty{0.7}{\percent}, whereas $c_\mathrm{out}$ decreased from \qty{45.38}{\mol\per\cubic\meter\per\second} to \qty{38.7}{\mol\per\cubic \meter \per \second}. However, the overall pressure drop of the packing increased by \qty{27.7}{\percent} from \qty{7.65}{\pascal} to \qty{9.77}{\pascal}.

\begin{Remark}
    The shape optimization approach presented above does only consider the mass transfer in the simplified simulation model and does not take into account the pressure drop over the packing, which is a crucial feature for practical applications. This could, e.g., be treated with a multi-criteria shape optimization, where the two conflicting goals of increasing the mass transfer and minimizing the pressure drop could be treated simultaneously, analogously to the approach considered in \cite{Blauth2024Multi}. Such an approach would try to approximate the Pareto-front, which consists of the optimal compromises between the conflicting targets. Here, an optimal compromise would be a packing design which cannot be further improved w.r.t.\ the mass transfer without increasing the pressure loss and vice versa. For a detailed introduction to multi-criteria optimization the reader is referred to, e.g., \cite{Ehrgott2005Multicriteria}. Future research efforts will be focused on this topic.
\end{Remark}

\subsection{Experimental Validation of the Optimized Design}
\label{ssec:experimental_results}

\begin{figure}[b]
    \centering
    \begin{subfigure}[b]{0.49\textwidth}
        \includegraphics[width=\textwidth]{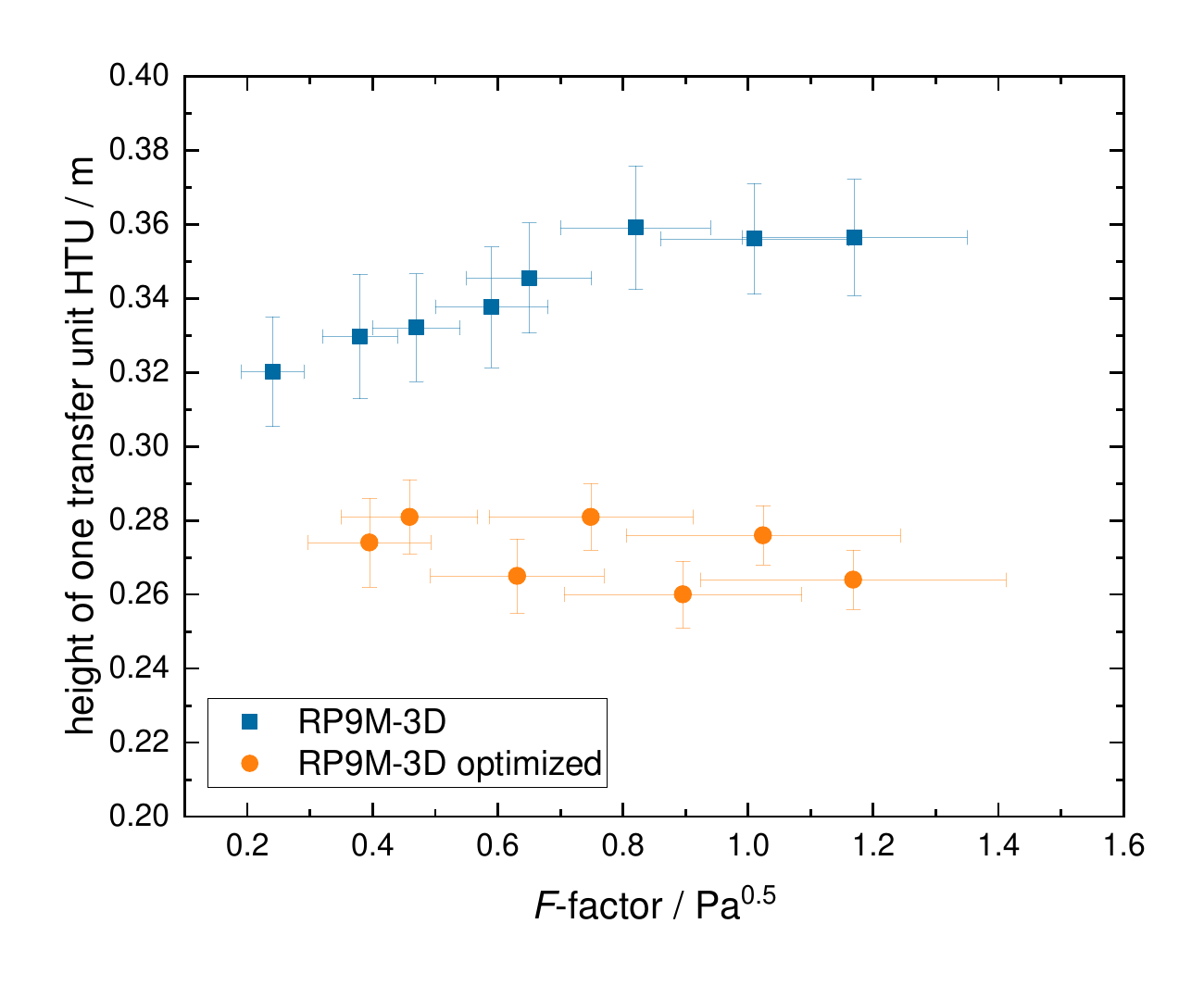}
        \caption{Separation efficiency.}
        \label{fig:experimental_HTU}    
    \end{subfigure}
    \begin{subfigure}[b]{0.49\textwidth}
        \includegraphics[width=\textwidth]{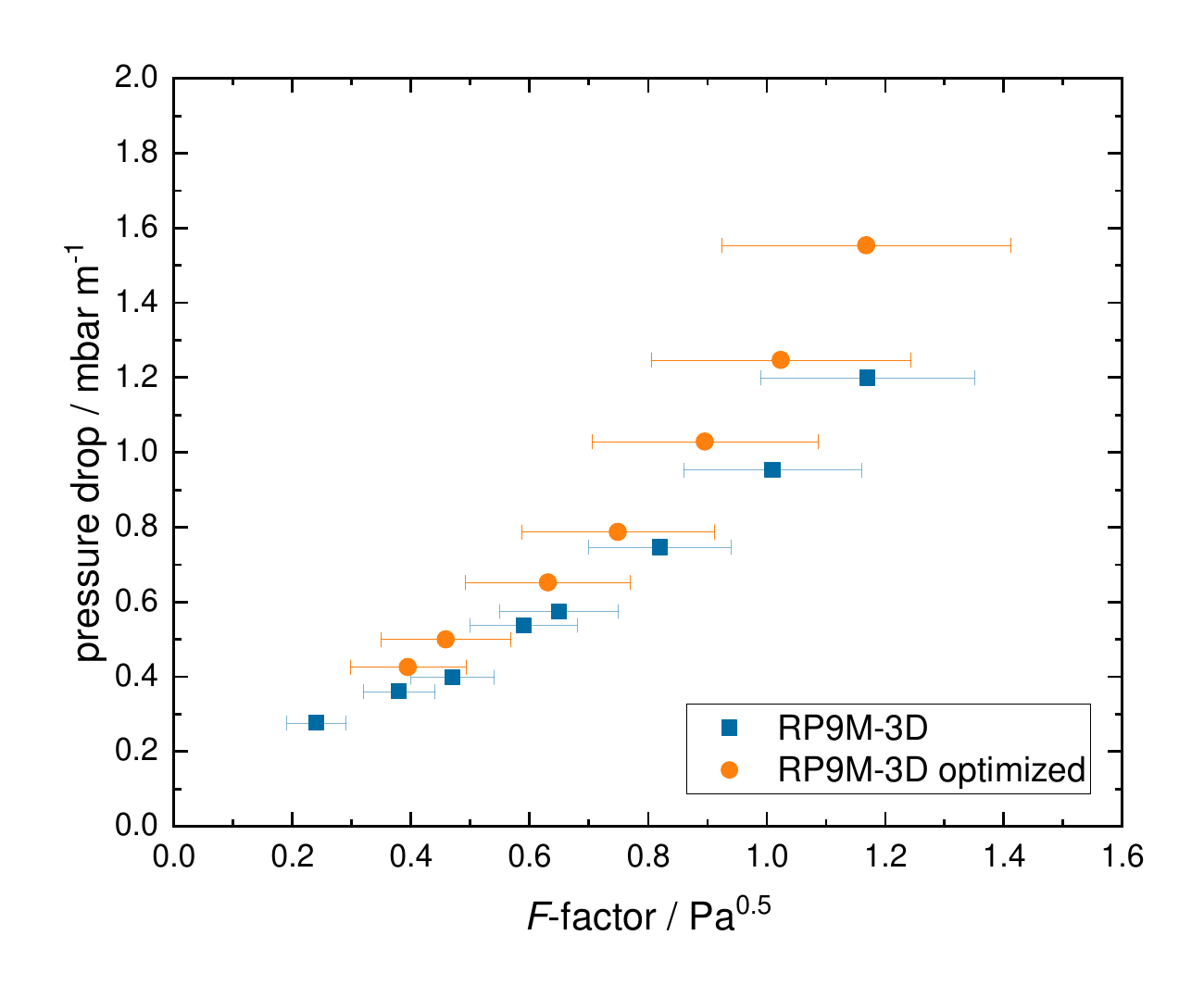}
        \caption{Wet pressure drop.}
        \label{fig:experimental_dP}
    \end{subfigure}
    \caption{Comparison of the experimentally measured separation efficiency (a) and wet pressure drop (b) of the original RP9M-3D and the optimized version.}
\end{figure}

To validate the simulation results the optimized packing was 3D-printed and the separation efficiency was measured in the laboratory-scale distillation column.
To ensure that the packings are as similar as possible, both the initial and optimized packings were manufactured by the same external additive manufacturing company using the same materials and post-processing steps.
The results and the comparison of the separation efficiency of the original and optimized RP9M-3D are shown in \cref{fig:experimental_HTU}. Additionally, the wet pressure drop is displayed in \cref{fig:experimental_dP}.

In \cref{fig:experimental_HTU} it can be observed that the height of one transfer unit (HTU) is significantly reduced for the optimized structure compared to the initial structure of the RP9M-3D which directly translates into higher separation efficiency. The improvement is around \qty{20}{\percent} at a \textit{F}-factor of $F=\qty{1}{\pascal^{0.5}}$. This is in good alignment with the simulation results presented in \cref{ssec:simulation_results}.

Regarding the pressure drop, which is presented in \cref{fig:experimental_dP}, the prediction from the simulation is met. Comparing the results again at $F=\qty{1}{\pascal^{0.5}}$, the experimentally determined wet pressure drop is approximately \qty{31}{\percent} higher for the optimized packing structure compared to the standard version. Even though this increase in pressure drop is significant, the absolute value of the pressure drop for the optimized RP9M-3D is still low enough to not limit the usual operation range of laboratory-scale distillation columns at $F<\qty{1}{\pascal^{0.5}}$ \cite{neukaufer_development_2022}.

\section{Conclusions and Outlook}
\label{sec:conclusions}

The shape optimization of a structured packing for a distillation column has been investigated with the aim of increasing its separation efficiency. After discussing the experimental setup and simulation methods, the optimized and initial designs have been compared. The results show that the predicted increase in separation efficiency using the simulation model is in very good agreement with the experimentally measured separation efficiency, yielding an increase in performance of about \qty{20}{\percent}. Hence, the presented approach of using shape optimization to improve the packing design seems to be extremely promising.

Applying and extending the presented shape optimization framework to other types of structured packings, which have been previously optimized manually \cite{neukaufer_development_2022}, to further increase their separation efficiency is an interesting direction for future work. Moreover, the shape optimization should also take into account the pressure loss in the packing with the goal of maximizing mass transfer while minimizing the pressure loss, which would naturally lead to multi-criteria shape optimization problems. Finally, the surrogate model for the distillation process could be further improved to better reflect the underlying physical and chemical effects without increasing the computational complexity of the problem too much, allowing for a more detailed shape optimization.

\bibliographystyle{elsarticle-num}

\appendix
\section{Details for the Simulations}
\label{app:solver_settings}

\subsection{Details for the Ansys\textsuperscript{\textregistered} Fluent Simulation}
\label{app:fluent_settings}

The computational mesh is generated using Ansys\textsuperscript{\textregistered} Meshing and consists of about \num{2.35}~million nodes and \num{13.7}~million tetrahedrons. The cell size is adapted based on the curvature of the geometry, so that the smallest mesh elements have an edge length of about \qty{6e-5}{\meter} and the average edge length of elements is about \qty{2.5e-4}{\meter}. 
For the simulation, the pressure-based solver with the default solver settings of Ansys\textsuperscript{\textregistered} Fluent was used. In particular, the coupled approach is used for solving the saddle-point problem of the incompressible Navier-Stokes equations \eqref{eq:navier_stokes} and second order upwind fluxes are used for all variables. Moreover, the convection-diffusion equation \eqref{eq:concentration} is implemented as a so-called user-defined scalar. Fluent uses a pseudo time stepping to compute the steady state solution of \eqref{eq:navier_stokes} and \eqref{eq:concentration}. For this, using \num{1000} pseudo time steps was sufficient for decreasing all residuals by nearly five orders of magnitude, indicating a good numerical convergence.

\subsection{Details for the FEniCS Simulation}
\label{app:fenics_settings}

As stated in Section~\ref{ssec:shape_optimization}, the incompressible Navier-Stokes equations \eqref{eq:navier_stokes} and the convection-diffusion equation \eqref{eq:concentration} are discretized with piecewise linear Lagrange elements for all variables. To stabilize the discretized equations, the following methods are employed: For the Navier-Stokes equations, a pressure-stabilized Petrov-Galerkin (PSPG) stabilization is used to circumvent the instability caused by the equal-order elements for the saddle-point structure of \eqref{eq:navier_stokes}. To stabilize the influence of convection in the momentum equation, a streamline-upwind Petrov-Galerkin (SUPG) stabilization is employed. For enhancing mass conservation due to the incompressibility, a grad-div or least-squares incompressibility constraint (LSIC) stabilization is used. For a detailed overview of these stabilization techniques for the Navier-Stokes equations, the reader is referred, e.g., to \cite{John2016Finite,Elman2014Finite}. 
For stabilizing the convection-dominated transport equation equation \eqref{eq:concentration}, both SUPG and crosswind stabilization are employed, where the formulation for both methods is taken from \cite{Codina1993discontinuity}. For more details regarding numerical stabilization of transport equations with the finite element method, the reader is referred, e.g., to \cite{Elman2014Finite,Roos1996Numerical}.

To solve the resulting discretized equations, an inexact Newton's method described in \cite{Kelley1995Iterative}, which is implemented in cashocs, is used with a relative tolerance of \num{1e-5}. The solution of the resulting linear systems is done with PETSc \cite{Balay2020PETSc}, where a GMRES (generalized minimal residual) method is used as iterative solver. For the Navier-Stokes equations \eqref{eq:navier_stokes} the so-called PCFIELDSPLIT preconditioner implemented in PETSc is employed, which is based on the Schur complement of the saddle-point problem. This corresponds to a SIMPLE-type (semi-implicit method for pressure-linked equations) preconditioner as discussed, e.g., in \cite{Elman2008taxonomy}. Both the momentum and the Schur complement block are preconditioned with BoomerAMG \cite{Henson2002BoomerAMG}, an algebraic multigrid method implemented in HYPRE \cite{hypre} and interfaced through PETSc. For the convection-diffusion equation \eqref{eq:concentration}, a GMRES method preconditioned by the algebraic multigrid method BoomerAMG is used.

\end{document}